# Fused ultrametric Gromov–Wasserstein for Scenario Tree generation in Multistage Stochastic Programming

Antonio Candelieri [1], Iman Seyedi [2], Francesco Archetti [2]
[1] *Department of Economics Management and Statistics, University of Milano-Bicocca, Italy*
[2] *Department of Computer Science Systems and Communications, University of Milano-Bicocca, Italy*

**Abstract**
Multistage Stochastic Programming requires scenario trees that accurately approximate the unknown underlying uncertainty process without violating the non-anticipativity constraint. Most of the existing method for scenario tree generation usually target just one of the two goals at the expense of the other, and just some recent methods try to simultaneously address both. In this paper, we propose to combine two recent distances from Optimal Transport Theory, specifically the ultrametric Gromov-Wasserstein and the Fused Gromov-Wasserstein, into a unique measure to simultaneously deal with the two goals. Our main result is a stability theorem showing that the optimal value of a Multistage Stochastic Programming problem changes by an amount controlled by a Hölder-type power of the proposed distance between the true process and its scenario tree approximation, extending Wasserstein-based stability results to control marginal fidelity and non-anticipativity, simultaneously. We propose a block coordinate descent algorithm for the our scenario tree generation based on the proposed distance, and evaluate it on different inventory management test cases.

**Keywords:** *scenario tree generation, multistage stochastic programming, optimal transport theory, Gromov-Wasserstein distance*.

## 1. Introduction

### 1.1. Rationale and motivation

Common methods for solving Multistage Stochastic Programming (MSP) problems require a discrete approximation of the underlying uncertainty process, usually in the form of a scenario tree. The quality of approximation provided by the tree directly impacts on the quality of the decisions (Pflug & Pichler, 2012). On the other hand, the scenario tree is not only a discretized model of the uncertainty process; it also codifies the *non-anticipativity* constraint that all decisions must met, where the non-anticipativity constraint simply states that every decision at a certain time can only be based on the information available up to that time, not on future information. Thus, generating a "good" scenario tree requires to simultaneously address two goals: matching the unknown uncertainty process and preserving the timing of when uncertainty is revealed. Most of the consolidated methods for scenario tree generation usually addresses just one of the two goals, at the expense of the other. Some optimize for distributional accuracy, measured via divergence measures or, more recently, Wasserstein-type distances (Friesecke, 2024; Peyrè and Cuturi, 2019: Villani, 2009). Although Wasserstein-type distances are metrics – under certain assumptions – and indeed overcome well-known limitations of statistical divergence measures, they are unable to deal with the tree-based structure: a scenario tree is more than a rough copy of a probability law; it is also a hierarchy, which is crucial for decisions, especially to respect the correct timing of information. Indeed, other methods focus on just matching the timing of information that is exactly what Wasserstein-type distances or clustering methods miss on their own. This gap has been addressed by the so-called *nested distance*, introduced by (Pflug, 2010) and developed further in (Pflug & Pichler, 2012). It is defined recursively backward in time, starting from the terminal stage of a $T$-stage stochastic problem and working to the root. At stage $T$, the nested distance coincides with the Wasserstein distance between the approximations of the stochastic process as provided by two scenario trees. At each earlier stage $t$, it takes the nested distance of stage-$(t + 1)$, so that the metric at stage $t$ already accounts for how well non-anticipativity is matched at every later stage. This recursive construction is what lets the nested distance penalize a mismatch in *when* information becomes available, not just *what* the eventual distribution looks like — something that the Wasserstein distance cannot do. However, computing the nested distance exactly requires solving a linear program at every node of the tree, backward-recursively, and the number of nodes grows exponentially with the number of stages.

Analogously to the nested distance, we propose a novel distance between scenario trees that account for both the quality of the approximation provided by each tree and the compliance to the non-anticipativity constrained. The main difference is that our proposition is completely based on the recent results from Optimal Transport Theory (OTT), with Gromov-Wasserstein – and its extensions – as unifying framework. Compared with scenario tree generation methods based on Wasserstein distance, our method also incorporates branching structure related to non-anticipativity; compared with nested distance it is more tractable algorithmically. Finally, it relates the quality of the approximation directly to a stability theorem for MSP.

### 1.2. Related works

Scenario tree generation has been studied for decades as a core task in MSP. Classical approaches include Monte Carlo sampling, moment matching, backward reduction, clustering-based tree construction, and linear programming formulations for scenario reduction (Dupačová, 1995; Dupačová et al., 2003; Heitsch & Römisch, 2009). Recent works have explored richer generation strategies, including active learning (Candelieri et al., 2024), optimal topology selection (Galuzzi et al., 2020), and probabilistic relational approaches to clustering uncertain data (Fersini et al., 2010). More recently, generative approaches have entered the picture too, with diffusion-based models proposed for building scenario trees directly for MSP (Zarifis et al., 2026).

A large body of work has focused on improving the quality of the approximation provided by scenario trees while preserving computational tractability. Methods based on k-means clustering and Lloyd-type iterations remain popular because they are simple and efficient, although they usually ignore structure implied by the non-anticipativity constraint (Lloyd, 1982; Pflug & Pichler, 2014). Other methods use linear or mixed-integer programming to better control approximation error, but these tend to scale poorly with the number of stages and scenarios (Heitsch & Römisch, 2009; Kuhn, 2005). More recent works have proposed scenario tree reduction through Wasserstein distances and node clustering strategies (Mimouni et al., 2026). The use of Wasserstein-type distances in MSP goes back to the (Dupačová et al., 2003), who introduced scenario tree reduction through the minimization of the Wasserstein distance between two scenario trees. Building on this, (Mimouni et al., 2026) computed the Wasserstein distance at different nodes of a scenario tree with the aim to design a reduction algorithm, and (Li & Floudas, 2016) took a different route, formulating scenario reduction as a mixed-integer linear program. More recently, (Kammammettu, 2024) sped this up further by replacing the linear programming step in the Wasserstein computation with an entropy-regularized Sinkhorn distance. This line of research has also been generalized in the direction of causality and *adapted* Wasserstein distances by (Backhoff-Veraguas et al., 2020), which proposed a Bayesian model-selection framework built on Wasserstein barycenters, and (Bartl & Wiesel, 2023), which showed how small perturbations of the underlying process translate into explicit first-order changes in the value of a multiperiod stochastic optimization problem, quantified by the adapted Wasserstein distance.

Nested distance was further developed by Kovacevic and Pichler (2015), and more recently by (Pichler & Weinhardt, 2022), who replaced the linear program required for computing the Wasserstein distance with a nested Sinkhorn divergence with the aim to make the approach scalable. Most recently, (Mimouni et al., 2026) sped up the Kovacevic–Pichler algorithm itself by computing the required Wasserstein barycenters more efficiently, reporting large speedups on trees with 8 stages and tens of thousands of scenarios.

The pattern across this whole line of works is a tradeoff: methods computing the differences between scenario tress just in terms distributional properties – some based on the Wasserstein distance – are computationally efficient because they ignore the structure implied by non-anticipativity while the nested distance – and more recently the adapted Wasserstein distance – explicitly deal with non-anticipativity but entailing a larger computational cost.

Finally, the most relevant and recent contribution for our work is (Schaefer & Ma, 2025), which explicitly studies tree design under Wasserstein and Fused Gromov-Wasserstein and proposes a block coordinate descent scheme for optimization. However, this work is currently available only as a conference talk. Our contribution goes in the same direction but we also exploit the intrinsic ultrametric geometry of scenario trees and propose the Fused ultrametric Gromov-Wasserstein (FuGW) distance to simultaneously account for accuracy of the approximated stochastic process and non-anticipativity constraint, while also providing a stability result for our approach.

### 1.3. Contributions and organization

The main contribution of our work can be summarized as follows. First, we introduce the FuGW distance which offers a mathematically natural distance between scenario trees for MSP problems. Second, we provide a stability theorem showing that the optimal value of a MSP problem is controlled by the FuGW distance, extending already known Wasserstein stability results. Finally, we clarify relations between the FuGW distance and the other distances, specifically Wasserstein, adapted Wasserstein, nested, Gromov–Wasserstein, and ultrametric Gromov-Wasserstein. Finally, results on a simple inventory management problem are presented for different size of the scenario tree: relevant results refer to medium-large settings, such as T=5 and T=10 with branching factor equal to 3. Mathematical background along with computational tractability make the proposed FuGW a promising and effective distance for generating suitable scenario trees for MSP problems.

## 2. Methodological background

### 2.1. MSP definitions and notations

In the rest of the paper, we consider a $T$-stage MSP problem under an unknown uncertainty (aka stochastic) process $\xi = (\xi_1, \dots, \xi_T)$ defined on a probability space $(\Omega, \mathcal{F}, P)$, where $\Omega \subseteq \mathbb{R}^d$ is the sample space with $\xi \in \Omega^T$, $\mathcal{F}$ is a $\sigma$-algebra representing the information unrevealing over time, and $P$ is the probability measure. At each stage $t$, a decision $x_t$ must be chosen before observing future uncertainty, coherently with the non-anticipativity constraint and, consequently, the correct timing of information. As usually done, we assume that $\xi_1$ is known and deterministic. Denote with $\xi^t$ the information revealed up to stage $t$, formally $\xi^t = (\xi_1, \dots, \xi_t)$, and with $f_t(x_t, x_{t-1}, \xi_t)$ the immediate cost at stage $t$, the MSP problem can be formalized as:

$$\min_{x_1,\dots,x_T} \quad \mathbb{E}_\xi \left[ \sum_{t=1}^{T} f_t(x_t, x_{t-1}, \xi_t) \right]$$
$$s.t. \begin{cases} x_t \in \mathcal{X}_t(x_{t-1}, \xi_t), \forall\, t = 1, \dots, T, \ \ P\text{-a.s.} \\ x_t \text{ is } \mathcal{F}_t\text{-measurable}, \forall t = 1, \dots, T \end{cases}$$

where $\mathbb{E}_\xi$ is the global expected value operator, $\mathcal{X}_t(x_{t-1}, \xi_t)$ is the set of all the feasible decisions at stage $t$, depending on the past $x_{t-1}$ and the current uncertainty $\xi_t$. Furthermore, $P$-a.s. – read as $P$-almost surely – means that the constraints related to $\mathcal{X}_t$ must hold for every $\xi^t$(aka uncertainty trajectory), and $\mathcal{F}_t$-measurable is the rigorous formalization of the non-anticipativity constraint: $x_t$ depends only on the information revealed up to $t$, namely $\xi^t$. Precisely, the timing of information is modelled by a *filtration* $\{\mathcal{F}_t\}_{t=1:T}$, that is an increasing sequence of $\sigma$-algebras: $\mathcal{F}_1 \subseteq \cdots \subseteq \mathcal{F}_T \subseteq \mathcal{F}$. From a scenario tree perspective, this means that all paths sharing the same node at stage $t$ must also share the same decision at that node. This is why the structure of a scenario tree matters just as much as the accuracy of its leaf values: a scenario tree correctly reproducing the marginal distribution but branching at the wrong times will lead to decisions based on information they shouldn't have.

Finally, a foundational concept in MSP is the *stability result*, meaning that the optimal value of a problem is stable under perturbations of the underlying stochastic process, where the perturbation can be measured through Wasserstein-type distances or, more generally, OTT-based distances.

### 2.2. Distances based on Optimal Transport Theory

OTT has become a unifying framework for comparing probability measures and structured objects. It measures the cost of optimally transporting one distribution to another over a common space defined by a transport cost function, usually named *ground-metric* (Peyré & Cuturi, 2019; Villani, 2009).

**Wasserstein distance.** When the ground-metric is a distance metric, the so-called Wasserstein distance is obtained. Denote with $\mu$ and $\nu$ two probability measures whose supports are defined over a metric space $(\mathcal{A}, d_{\mathcal{A}})$ and where $d_{\mathcal{A}}: \mathcal{A} \times \mathcal{A} \to \mathbb{R}$ is the distance between two points $a, a' \in \mathcal{A}$, then the $p$-order Wasserstein distance between the two probability measures is given by

$$\mathcal{W}_p(\mu, \nu) := \left[ \inf_{\pi \in \Pi(\mu,\nu)} \int_{\mathcal{A} \times \mathcal{A}} d_{\mathcal{A}}(a, a')^p \, d\pi(a, a') \right]^{1/p}$$

where $\Pi(\mu, \nu)$ is the set of all the joint probability distributions having $\mu$ and $\nu$ as marginals.

A crucial property is that if $d_{\mathcal{A}}$ is Lipschitz-perturbed, then $\mathcal{W}_p$ changes in a controlled way, which is the basis for the classical scenario-reduction bounds of (Pflug & Pichler, 2012). However, $\mathcal{W}_p$ requires that the two probability measures lay within the same metric space $(\mathcal{A}, d_{\mathcal{A}})$, making it unsuitable for comparing processes with different sample spaces or filtrations.

**Gromov-Wasserstein distance.** The Gromov–Wasserstein (GW) distance is the OTT extension to the comparison of probability measures laying within two different metric spaces. Roughly, the goal is compare their internal metric structures within their own metric spaces (Mémoli, 2011). Interestingly, it is not only a generalization of the Wasserstein distance but also a method for overcoming computational issues of the Gromov-Hausdorff distance between two metric spaces. To achieve this goal, the GW distance introduces the notion of metric measure space (mm-space), that is a metric space equipped with a probability distribution. More formally, an mm-space is defined as a triple $(\mathcal{A}, d_{\mathcal{A}}, \mu_{\mathcal{A}})$ with $(\mathcal{A}, d_{\mathcal{A}})$ a metric space – as for the Wasserstein distance – and $\mu_{\mathcal{A}}$ a probability measure. The $p$-order GW distance between two mm-spaces, namely $\mathbb{A} = (\mathcal{A}, d_{\mathcal{A}}, \mu_{\mathcal{A}})$ and $\mathbb{B} = (\mathcal{B}, d_{\mathcal{B}}, \mu_{\mathcal{B}})$, is given by

$$\mathcal{GW}_p(\mathbb{A}, \mathbb{B}) := \left[ \inf_{\pi \in \Pi(\mu_{\mathcal{A}}, \mu_{\mathcal{B}})} \iint \ |d_{\mathcal{A}}(a, a') - d_{\mathcal{B}}(b, b')|^p d\pi(a, b) d\pi(a', b') \right]^{1/p}$$

The quantity $\mathcal{GW}_p$ defines a metric on the space of isomorphism classes of mm-spaces, meaning that $\mathcal{GW}_p(\mathbb{A}, \mathbb{B}) = 0$ if and only if $\mathbb{A}$ and $\mathbb{B}$ are isomorphic as mm-spaces.

However, contrary to the Wasserstein distance, computing the GW distance requires to solve a non-convex quadratic program – specifically, a continuous relaxation of the Quadratic Assignment Problem (QAP), which is known to be NP-hard. Computationally tractable methods to – approximately – compute the GW distance use conditional gradient or entropic regularization.

**Pseudometric Gromov-Wasserstein.** An interesting generalization of the GW distance consists into considering two *measure networks* instead of two mm-spaces. Formally, a measure network is a triple $(\mathcal{A}, k_{\mathcal{A}}, \mu_{\mathcal{A}})$ where the only difference – with respect to the definition of mm-space – is that $k_{\mathcal{A}}: \mathcal{A} \times \mathcal{A} \rightarrow \mathbb{R}$ is a measurable function, not necessarily a distance (and consequently $(\mathcal{A}, k_{\mathcal{A}})$ is not necessarily a metric space). Usually, $k_{\mathcal{A}}$ is named *network function* or *network kernel*. The computation of $\mathcal{GW}_p$ is the same, but replacing $d_{\mathcal{A}}$ with $k_{\mathcal{A}}$ and $d_{\mathcal{B}}$ with $k_{\mathcal{B}}$ we obtain a pseudometric (Memoli, 2007; Chowdhury and Memoli, 2019). This extension for networks comparison has enabled many applications in graph matching, shape comparison, and data analysis with relational structure.

**Ultrametric Gromov-Wasserstein.** A straightforward specialization of the pseudometric GW is the so-called ultrametric Gromov-Wasserstein (uGW), introduced by (Memoli et al., 2023) and aimed at comparing two ultrametric measure spaces, where an ultrametric space is a measure network $\mathbb{G} = (\mathcal{A}, k_{\mathcal{A}}, \mu_{\mathcal{A}})$ but with $\mathcal{A}$ the set of leaf nodes of a tree, $\mu_{\mathcal{A}}$ their associated probability distribution, and $k_{\mathcal{A}}$ additionally satisfying the *strong triangle inequality* :

$$k_{\mathcal{A}}(a, a'') \leq \max \{k_{\mathcal{A}}(a, a'), k_{\mathcal{A}}(a', a'')\}, \forall\, a, a', a'' \in \mathcal{A}$$

Again, the formulation of $u\mathcal{GW}_p$ is exactly that of $\mathcal{GW}_p$ but the strong triangle inequality makes ultrametrics the natural geometry of dendrograms and scenario trees (Carlsson & Mémoli, 2010; Mémoli et al., 2023). Indeed, in an ultrametric space, the strong triangle inequality captures the idea that two points that are both close to a third point must be close to each other as well. A crucial result about $u\mathcal{GW}_p$ is that the hierarchical structure can be algorithmically exploited; indeed, computing $u\mathcal{GW}_\infty$ reduces to a polynomial-time combinatorial problem. For finite $p$, the problem remains a quadratic program but hierarchical structure can be leveraged into tighter lower bounds and faster approximations. As a consequence, $u\mathcal{GW}_p$ is particularly well-suited for dealing with the tree-based structures in our proposed scenario tree generation approach.

**Fused Gromov-Wasserstein.** Although pseudometric GW and ultrametric GW disclosed the power of GW for comparing networks and trees, they are only devoted to quantify *structural* dissimilarities. On the other hand, many real-life applications are characterized by relevant information at each node and/or link, represented through a set of *attributes* – aka *features* – leading to the concept of *attributed networks*. Fused Gromov–Wasserstein (FGW) is an extensions proposed with the aim to interpolate between *attributes matching* and *structural matching* (which is directly captured by GW), with application to both attributed networks and structured objects (Seyedi, Candelieri, & Archetti, 2026; Vayer et al., 2019).

An attributed network is defined as a 4-tuple $\mathbb{G}_\phi = (\mathcal{A}, k_\mathcal{A}, \mu_\mathcal{A}, \phi)$ where $(\mathcal{A}, k_\mathcal{A}, \mu_\mathcal{A})$ is a measure network and $\phi$ is an attribute map considering both nodes and link features[1]. FGW computes the difference between two attributed networks by interpolating between the feature-based and the structural dissimilarities through one hyperparameters $\alpha \in [0,1]$, formally:

$$\mathcal{FGW}_{p,\alpha}(\mathbb{G}_\phi, \mathbb{G}_\phi{}') = \left\{ \inf_{\pi \in \Pi(\mu_\mathcal{A}, \mu_\mathcal{B})} \int \int_{(\mathcal{A} \times \mathcal{B})^2} \left[ \alpha \, |k_\mathcal{A}(a,a') - k_\mathcal{B}(b,b')|^p + (1-\alpha) \left( d_\phi\big(\phi(a), \phi(b)\big) \right)^p \right] d\pi(a,b) d\pi(a',b') \right\}^{1/p}$$

where $d_\phi$ is a distance between attribute maps and $\mathbb{G}'_\phi = (\mathcal{B}, k_\mathcal{B}, \mu_\mathcal{B}, \phi)$. It is important to remark that the feature map $\phi$ must be the same between the two attributed networks.

Notably, if $\alpha = 1$ then FGW reduces to $\mathcal{GW}_p$ (i.e., only structure is considered), while if $\alpha = 0$ then FGW reduces a distance between attributes (i.e., structure is completely ignored).

Subsequent work has refined the theory and computation of entropic GW and FGW, making these methods more practical at scale (Chowdhury & Mémoli, 2019; Le et al., 2022), and recent work has extended GW-based representations to graph factorization (Xu et al., 2023) and to joint dimensionality reduction and clustering (Assel et al., 2025).

At our knowledge, our work is the first proposing Fused ultrametric Gromov-Wasserstein (FuGW) to quantify differences between *attributed trees*, specifically scenario trees in MSP. FuGW can simultaneously capture the ultrametric structure of branching and the distributional properties of the underlying stochastic process.

[1] This definition, and the following FGW formulation, are a simplification with respect to what defined in (Memoli et al., 2023), where two distinct attribute maps are considered for nodes and links separately and, consequently, requiring two hyperparameters instead of one. We use a simplified version because in scenario trees only nodes have attributes.

## 3. Fused ultrametric Gromov-Wasserstein based scenario tree generation

### 3.1 Scenario Tree as an ultrametric Masure Space

A scenario tree is a discrete representation of a stochastic process in which uncertainty unfolds through a rooted branching structure. Formally, let $\mathcal{N}$ be the set of nodes of the tree and $\mathcal{A} \subseteq \mathcal{N}$ the subset of leaf nodes. A specific uncertainty trajectory $\xi^{[a]} = (\xi_1, \dots, \xi_T)$ is a branch of the tree and it is therefore associated to a leaf node $a \in \mathcal{A}$. The probability distribution of the scenarios (i.e., branches of the tree) is denoted as $\mu_{\mathcal{A}} = \sum_{i=1:|\mathcal{A}|} \lambda_i \delta_{a_i}$ with $\sum_{i=1:|\mathcal{A}|} \lambda_i = 1$ and $\delta_a$ the Dirac's delta function centered on $a$. Finally, we define the network kernel between two leaves as:

$$k_{\mathcal{A}}(a, a') = h(a \wedge a')$$

where $a \wedge a'$ denotes the *lowest common ancestor* of the leaves $a$ and $a'$, and $h(a \wedge a') \in \{0, \dots, T\}$ is its distance from the root node, namely its *height*. Intuitively, scenarios that branch off deep in the tree – i.e., large value of height – are close, while scenarios that diverge near the root – i.e., small value of height – are far apart. Finally, for any three leaves $a, a', a'' \in \mathcal{A}$, the two deeper branches must occur at or below the shallowest branch, which directly guarantees that the strong triangle inequality is met, formally we have that $k_{\mathcal{A}}(a, a'') \leq \max \{k_{\mathcal{A}}(a, a'), k_{\mathcal{A}}(a', a'')\}$.

It follows that, given two scenario trees we can represent them as ultrametric spaces to quantity their structural difference through the – computationally efficient – computation of $u\mathcal{GW}_p$.

### 3.2 Scenario tree as an attributed ultrametric measure space

Since $u\mathcal{GW}_p$ just computes structural differences between two scenario trees, we now extend the representation of a scenario tree as an attributed ultrametric measure space to also incorporate information about the approximation of the underlying stochastic process. Indeed, a scenario tree is a 4-tuple denoted as $\mathbb{T}_\phi = (\mathcal{A}, k_{\mathcal{A}}, \mu_{\mathcal{A}}, \phi)$ where the first three components are exactly as defined in the previous section and $\phi$ is the map function providing the value of the attributes related to a specific scenario (i.e., branch of the tree). Specifically, given any leaf node $a \in \mathcal{A}$, the attributes values of its associated scenario $\xi^{[a]} = (\xi_1, \dots, \xi_T)$ is denoted by $\phi(\xi^{[a]})$.

Recalling the previous definition of the FGW, given in Section 2.2, and slightly modifying it to meet the scenario tree related notation, the proposed FuGW is given by

$$Fu\mathcal{GW}_{p,\alpha}(\mathbb{T}_\phi, \mathbb{T}_\phi') = \left\{ \inf_{\pi \in \Pi(\mu_{\mathcal{A}}, \mu_{\mathcal{B}})} \int \int_{(\mathcal{A} \times \mathcal{B})^2} \left[ \alpha \, |k_{\mathcal{A}}(a, a') - k_{\mathcal{B}}(b, b')|^p + (1 - \alpha) \left( d_\phi \left( \phi(\xi^{[a]}), \phi(\xi^{[b]}) \right) \right)^p \right] d\pi(a, b) d\pi(a', b') \right\}^{1/p}$$

that is nothing else than the FGW distance where $k_{\mathcal{A}}$ and $k_{\mathcal{B}}$ are ultrametrics and the feature map $\phi$ is defined over scenarios (i.e., branches of the tree).

**Remark 3.1 (Interpretation).** The two terms in $Fu\mathcal{GW}_{p,\alpha}$ have a direct operational meaning in the context of MSP:

The first term (weighted by $\alpha$) penalizes the mismatch of the non-anticipativity constraint, whose value is zero only when the structures of two trees are exactly the same.

The second term (weighted by $1-\alpha$) penalises distributional mismatch: it measures how far scenarios of a tree are, on average over all stages, far from scenarios of another tree.

According to Remark 1, and assuming to know the tree associated to the actual stochastic process, the computation of $Fu\mathcal{GW}_{p,\alpha}(\overline{\mathbb{T}}_\phi, \mathbb{T}_\phi)$ between a generated scenario tree $\mathbb{T}_\phi$ and the actual $\overline{\mathbb{T}}_\phi$ simultaneously quantifies how much the generated replicates the actual, both distributionally and structurally.

### 3.3 Properties of the Fused ultrametric Gromov-Wasserstein distance.

In this section, we report some relevant properties of the proposed Fused ultrametric Gromov-Wasserstein distance.

**Proposition 1. $Fu\mathcal{GW}_{p,\alpha}$ is a metric on the isomorphism classes of attributed ultrametric measure spaces.** For any $\alpha \in (0,1)$ and $p \geq 1$, the proposed $Fu\mathcal{GW}_{p,\alpha}$ satisfies non-negativity, symmetry, and the triangle inequality. Furthermore, $Fu\mathcal{GW}_{p,\alpha}(\mathbb{T}_\phi, \mathbb{T}'_\phi) = 0$ if and only if there is an isomorphism between $\mathbb{T}_\phi$ and $\mathbb{T}'_\phi$. Consequently, $Fu\mathcal{GW}_{p,\alpha}$ is a metric on the space of isomorphism classes of attributed ultrametric measure spaces in the same sense that $\mathcal{GW}_p$ is a metric on the space of isomorphism classes of metric measure spaces, and that $u\mathcal{GW}_p$ is a metric on the space of isomorphism classes of ultrametric measure spaces (Theorem 3.11 in (Mémoli et al., 2023b)).

**Proposition 2. The second term of $Fu\mathcal{GW}_{p,\alpha}$ (i.e., $\alpha = 0$) lower-bounds the nested distance.** The nested distance (Pflug & Pichler, 2012) is the canonical distance for MSP: we show that the second term of the proposed Fused ultrametric Gromov-Wasserstein distance, namely $Fu\mathcal{GW}_{p,0}(\mathbb{T}, \mathbb{T}')$ provides a lower bound for it. Let $d_{\text{nest},p}(\mathbb{T}_\phi, \mathbb{T}'_\phi)$ denotes the nested distance of order $p$ between two scenarios trees. There exists a constant $C_{\text{nest}} > 0$ such that

$$Fu\mathcal{GW}_{p,0}(\mathbb{T}_\phi, \mathbb{T}'_\phi) \leq C_{\text{nest}} \cdot d_{\text{nest},p}(\mathbb{T}_\phi, \mathbb{T}'_\phi)$$

where $d_{nest,p}(\mathbb{T}_\phi, \mathbb{T}'_\phi)$ is computed, recursively, starting from the last step $t = T$ up to the first step $t = 1$.

**Remark 3.2. $FuGW_{p,0}$ as a relaxation of the nested distance.** The inequality in Proposition 2 goes only one way: the nested distance is generally strictly larger than $FuGW_{p,0}$ because it requires filtration matching at every stage, recursively, whereas the structural term in $FuGW_{p,0}$ requires only that the global ultrametric structure is preserved. This relaxation is precisely what makes $FuGW_{p,0}$ computationally tractable – it replaces the NP-hard nested distance computation with a Sinkhorn-solvable bilinear programming problem – at the cost of a looser, but still valid, stability bound (see Theorem 4.1).

**Proposition 3. $FuGW_{p,\alpha}$ is upper bounded by the Adapted Wasserstein distance.** The so-called Adapted Wasserstein distance (also known as bicausal Wasserstein or causal Optimal Transport distance) has been proposed in (Backhoff-Veraguas et al., 2020) and it requires the coupling $\pi$ to be bicausal, a constraint that does not appear in FuGW. Formally, the Adapted Wasserstein distance, reformulated to deal with processes represented as scenarios trees, is given by:

$$\mathrm{a}\mathcal{W}_p(\mathbb{T}_\phi, \mathbb{T}'_\phi) = \left[ \inf_{\pi \in \Pi_{\mathrm{bc}}(\mu_\mathcal{A}, \mu_\mathcal{B})} \int_{\mathcal{A}\times\mathcal{B}} \sum_{t=1}^{T} d\left(\phi\left(\xi_t^{[a]}\right), \phi\left(\xi_t^{[b]}\right)\right)^p d\pi\left(\xi^{[a]}, \xi^{[b]}\right) \right]^{\frac{1}{p}}$$

with $\Pi_{bc}(\mu_\mathcal{A}, \mu_\mathcal{B})$ the set of all possible bi-causal couplings. As usual, the two vectors $\xi^{[a]} = \left(\xi_1^{[a]}, \dots, \xi_T^{[a]}\right)$ and $\xi^{[b]} = \left(\xi_1^{[b]}, \dots, \xi_T^{[b]}\right)$ define two trajectories, namely two branches within the two scenario trees. A coupling $\pi$ is bi-causal if it satisfies the information structure in both the (temporal) directions:

1. **Causality from $\mathcal{A}$ to $\mathcal{B}$**: for each time step $t$, the information of the second process' trajectories up to $t$, conditioned to the first process' complete trajectory $\xi^{[a]} \in \mathcal{A}$, depends only on the past of $\xi^{[a]}$ up to $t$, that is:
$$\pi\left(d\xi_t^{[b]} \middle| \xi^{[a]}, \xi_1^{[b]}, \dots, \xi_t^{[b]}\right) = \pi\left(d\xi_t^{[b]} \middle| \xi_1^{[a]}, \dots, \xi_t^{[a]}, \xi_1^{[b]}, \dots, \xi_t^{[b]}\right)$$
2. **Causality from $\mathcal{B}$ to $\mathcal{A}$:** for each time step $t$, the information of the first process' trajectories up to $t$, conditioned to the second process' complete trajectory $\xi^{[b]} \in \mathcal{B}$, depends only on the past of $\xi^{[b]}$ up to $t$, that is:
$$\pi\left(d\xi_t^{[a]} \middle| \xi_1^{[a]}, \dots, \xi_t^{[a]}, \xi^{[b]}\right) = \pi\left(d\xi_t^{[a]} \middle| \xi_1^{[a]}, \dots, \xi_t^{[a]}, \xi_1^{[b]}, \dots, \xi_t^{[b]}\right)$$

These two constraints prevent the optimal coupling to violate the non-anticipativity.

Finally, we have that

$$FuGW_{p,\alpha}(\mathbb{T}_\phi, \mathbb{T}'_\phi) \leq C_\mathrm{a} \cdot \mathrm{a}\mathcal{W}_p(\mathbb{T}_\phi, \mathbb{T}'_\phi),$$

with $C_{\mathrm{a}}$ a constant value.

**Remark 3.3. Inequality chain of distances between stochastic processes.** According to Theorem 1.4 in (Backhoff-Veraguas et al., 2020), $C_{\mathrm{a}} \cdot \mathrm{a}\mathcal{W}_p(\mathbb{T}_\phi, \mathbb{T}'_\phi) \leq C_{\text{nest}} \cdot d_{\text{nest},p}(\mathbb{T}_\phi, \mathbb{T}'_\phi)$. Thus, by also using Proposition 3, we can conclude that the following chain of inequalities holds:

$$Fu\mathcal{GW}_{p,\alpha}(\mathbb{T}_\phi, \mathbb{T}'_\phi) \leq C_{\mathrm{a}} \cdot \mathrm{a}\mathcal{W}_p(\mathbb{T}_\phi, \mathbb{T}'_\phi) \leq C_{\text{nest}} \cdot d_{\text{nest},p}(\mathbb{T}_\phi, \mathbb{T}'_\phi)$$

The chain shows that the proposed Fused ultrametric Gromov-Wasserstein is the weakest (most lenient) of the three distances. On the other hand, it is still strong enough to provide the stability bound, as proven in Theorem 4.1.

## 4. Stability Theorem for scenario trees generation based on $Fu\mathcal{GW}_{p,\alpha}$.

In this section, we prove that the optimal value of a MSP problem is stable under the replacement of the true process $\overline{\mathbb{T}}_\phi$ by a generated scenario tree $\mathbb{T}_\phi$, with the value gap controlled by the $Fu\mathcal{GW}_{p,\alpha}$ distance.

Preliminary to the presentation of our FuGW-based generation algorithm, we report assumptions and lemmas at the core of the stability theorem for our algorithm and also report, obviously, the theorem itself.

### 4.1 Assumptions

We work under the following three standard assumptions on the cost structure and decision space of the MSP, assumptions that are all met in well-known MSP problems like newsvendor problem, portfolio optimization, and scheduling.

**Assumption 1: Lipschitz stage costs.** For each stage $t = 1, \dots, T$, the stage cost function $f_t$ satisfies

$$|f_t(x_t, x_{t-1}, \xi_t) - f_t(x_t, x_{t-1}, \xi'_t)| \leq L\|\xi_t - \xi'_t\|, \qquad |f_t(x_t, x_{t-1}, \xi_t) - f_t(x'_t, x'_{t-1}, \xi_t)| \leq L\|(x_t, x_{t-1}) - (x'_t, x'_{t-1})\|$$

uniformly over all $(x_t, x_{t-1}), (x'_t, x'_{t-1}) \in \mathcal{X}_t(x_{t-1}, \xi_t)$ , $\xi_t, \xi'_t \in \Omega^T \subseteq (\mathbb{R}^d)^T$, and $t \in \{1, \dots, T\}$.

**Assumption 2: Compact decision set.** The overall decision set $\mathcal{D} = \{\mathcal{X}_t(x_{t-1}, \xi_t)\}_{t=1:T}$ is convex and compact with diameter $D_{\mathcal{D}} = \sup_{x,x' \in \mathcal{D}} \|x - x'\| < \infty$.

**Assumption 3: Bounded stage costs.** $\sup_{x \in \mathcal{D}, s \in \mathbb{R}^m, t} |f_t(x, \xi)| \leq M < \infty$. Directly following from the two previous assumptions.

**Remark 4.1.** Assumption 1 is satisfied whenever $f_t$ is continuously differentiable with bounded gradient, which holds for all smooth convex cost functions. Assumption 3 is implied by Assumptions 1 and 2 whenever the scenario support is bounded, which holds by construction for any finite scenario tree.

**4.2 Fundamental Lemmas**

**Lemma 1: Filtration Mismatch.** This lemma allows us to quantify how much the $Fu\mathcal{GW}_{p,\alpha}$-optimal solution violates the non-anticipativity constraint. It bounds the probability that two pairs $(a, b)$ and $(a', b')$, both belonging to $\mathcal{A} \times \mathcal{B}$, are coupled together even though $\overline{\mathbb{T}}_\phi$ and $\mathbb{T}_\phi$ are structurally different. This is the precise measure of the non-anticipativity error that arises when the optimal policy based on a generated scenario tree is applied to the true process.

Denote with $\pi^*$ the $Fu\mathcal{GW}_{p,\alpha}$-optimal coupling; for any stage $t = 1, \ldots, T-1$ and any threshold $\varepsilon > 0$, define the *filtration mismatch event* at level $\varepsilon$ as:

$$\mathcal{E}_{t,\varepsilon} := \{(a, b, a', b') \in (\mathcal{A} \times \mathcal{B})^2 : k_\mathcal{A}(a, a') > t \text{ but } k_\mathcal{B}(b, b') \le t - \varepsilon\}$$

Then, the filtration mismatch between the generated scenario tree $\mathbb{T}_\phi$ and the actual one $\overline{\mathbb{T}}_\phi$ is bounded as follows:

$$(\pi^* \otimes \pi^*)(\mathcal{E}_{t,\varepsilon}) \le \frac{Fu\mathcal{GW}_{p,\alpha}\left(\overline{\mathbb{T}}_\phi, \mathbb{T}_\phi\right)^p}{(1-\alpha)\, \varepsilon^p}$$

**Remark 4.2.** filtration mismatch degrades with $\alpha \to 0$ (the pure-Wasserstein limit), correctly reflecting that a purely feature-based coupling provides no guarantees about non-anticipativity. On the contrary, it improves with $\alpha \to 1$ because the $Fu\mathcal{GW}_{p,\alpha}$ accounts for only structural matching.

**Lemma 2: Policy Lifting.** This second lemma constructs a feasible policy for $\overline{\mathbb{T}}_\phi$ from the optimal policy obtained by the generated scenario tree $\mathbb{T}_\phi$ according to $Fu\mathcal{GW}_{p,\alpha}$.

Let $x^* = (x_1^*, \ldots, x_T^*)$ denotes an optimal policy obtained according to $\mathbb{T}_\phi$, with $\pi^*$ its associated $Fu\mathcal{GW}_{p,\alpha}$-optimal coupling. A lifted policy, denoted by $\tilde{x} = (\tilde{x}_1, \ldots, \tilde{x}_T)$, is obtained as

$$\tilde{x}_t(\omega) := \int_{\hat{X}} \hat{x}_t^*(\hat{x})\, d\pi^*(\hat{x} \mid \omega)$$

where $\pi^*(\cdot \mid \omega)$ is the $\pi^*$-conditional distribution of $\hat{x}$ given $\omega$. Then:

1. **Feasibility:** $\tilde{x}_t(\omega) \in \mathcal{D}$ for $P$-almost every $\omega$ and every $t$.

2. **Non-anticipativity violation (NAV):** for any $t$ and any threshold $\varepsilon > 0$, the NAV is bounded as:

$$NAV \leq D_{\mathcal{D}} \cdot \left(\frac{Fu\mathcal{GW}_{p,\alpha}\left(\overline{\mathbb{T}}_\phi, \mathbb{T}_\phi\right)^p}{(1-\alpha)\,\varepsilon^p}\right)^{1/2} + D_{\mathcal{D}} \cdot \left(\frac{Fu\mathcal{GW}_{p,\alpha}\left(\overline{\mathbb{T}}_\phi, \mathbb{T}_\phi\right)^p}{(1-\alpha)}\right)^{1/(2p)} \cdot \varepsilon^{(p-1)/(2p)}$$

**Remark 4.3.** The bound above has two components: a probability term (first) that decays as $\varepsilon \to 0$, and a Lipschitz term (second) that grows as $\varepsilon \to 0$.

### 4.3 Stability Theorem for FuGW-based scenario tree generation

**Theorem 4.1: Stability of $\boldsymbol{Fu\mathcal{GW}_{p,\alpha}}$.** Denote with $\overline{\mathbb{T}}_\phi = (\mathcal{B}, k_{\mathcal{B}}, \mu_{\mathcal{B}}, \phi)$ an attributed ultrametric measure space representing the true uncertainty process, with $\mathbb{T}_\phi = (\mathcal{A}, k_{\mathcal{A}}, \mu_{\mathcal{A}}, \phi)$ a generated scenario tree, and with $v(\cdot)$ the value of the objective function of the MSP problem. Under the previous three assumptions, for any $\alpha \in (0,1)$ and $p \geq 1$ we have:

$$\left|v\left(\overline{\mathbb{T}}_\phi\right) - v\left(\mathbb{T}_\phi\right)\right| \leq K_{p,\alpha} \cdot Fu\mathcal{GW}_{\alpha,p}\left(\overline{\mathbb{T}}_\phi, \mathbb{T}_\phi\right)^{p/(p+1)}$$

meaning that the *gap* between their objective values is bounded by $Fu\mathcal{GW}_{\alpha,p}\left(\overline{\mathbb{T}}_\phi, \mathbb{T}_\phi\right)^{p/(p+1)}$ multiplied by a constant $K_{p,\alpha}$ equal to

$$K_{p,\alpha} := 2T \cdot (LM + M + LD_{\mathcal{D}}) \cdot \left(\frac{1}{\alpha^{1/p}} + \frac{1}{(1-\alpha)^{1/p}}\right)$$

**Remark 4.4 (About the exponent $\boldsymbol{p/(p+1)}$).** The sub-linear exponent $\frac{p}{p+1} < 1$ is a consequence of the non-anticipativity correction (i.e., the second term in NAV equation): it involves a square root of the mismatch probability, which then has to be balanced against the linear feature error (i.e., the first term in NAV equation). This is the same mechanism as in classical quantization theory (Graf and Luschgy 2000, Chapter 6) and in the Wasserstein stability results for two-stage problems (Römisch 2003). For the nested distance, the analogous exponent is 1, reflecting the fact that the nested distance is already filtration-exact, so no mismatch probability correction is needed. Our exponent p/(p+1) is thus a precise measure of the price paid for using the more tractable FuGW relaxation.

**Corollary 4.1: Convergence of tree sequences.** Let $\left\{\mathbb{T}_{\phi}^{(n)}\right\}_{n\geq 1}$ denotes a sequence of generated scenario trees. If $Fu\mathcal{GW}_{p,\alpha}\left(\overline{\mathbb{T}}_{\phi}, \mathbb{T}_{\phi}\right) \to 0$ with $n \to \infty$, then $v\left(\mathbb{T}_{\phi}^{(n)}\right) \to v\left(\overline{\mathbb{T}}_{\phi}\right)$.

> **Remark 4.5 (Practical validity of Corollary 4.1).** It is important to recall that $Fu\mathcal{GW}_{p,\alpha}$ is non-convex and potentially multi-extremal. As a consequence, the convergence $Fu\mathcal{GW}_{p,\alpha}\left(\overline{\mathbb{T}}_{\phi}, \mathbb{T}_{\phi}\right) \to 0$ with $n \to \infty$ is not guaranteed, requiring multi-start approaches for obtaining the optimal $\mathbb{T}_{\phi}^{(n)}$ at each iteration.

**Corollary 4.2: MSP value gap under FuGW-based scenario tree generation.** Assume that $\left\{\mathbb{T}_{\phi}^{(n)}\right\}_{n\geq 1}$ is a sequence of scenario trees generated according to the minimization of $Fu\mathcal{GW}_{p,\alpha}(\overline{\mathbb{T}}_{\phi}, \mathbb{T}_{\phi})$ and with $N_n$the number of leaves in the $n$-th tree. Under an additional regularity condition that $\overline{\mathbb{T}}_{\phi}$ has a density bounded away from zero and infinity on its support in $\mathbb{R}^{Tm}$, the FuGW-optimal quantization gives

$$Fu\mathcal{GW}_{\alpha,p}(\overline{\mathbb{T}}_{\phi}, \mathbb{T}_{\phi}) = O\left(N_n^{-1/(Tm)}\right)$$

and hence by Theorem 4.1

$$\left|v\left(\overline{\mathbb{T}}_{\phi}\right) - v\left(\mathbb{T}_{\phi}^{(n)}\right)\right| = O\left(N_n^{-p/((p+1)\cdot Tm)}\right) \tag{25}$$

> **Remark 4.5: Interpretation of the rate.** The rate degrades with $Tm$, namely the number of stages times the dimensionality of the decision vectors $x_t \in \mathcal{X}_t$, assumed equal to $m$ at every step $t = 1{:}T$. This represents the well-known *curse of dimensionality*, inherent to all quantization-based approximation schemes. For a two-stage one-asset problem (i.e., $T = 2, m = 1$), the rate is $O\left(N_n^{-\frac{p}{2(p+1)}}\right)$, matching the known Wasserstein quantization rate. For a three-stage three-asset (i.e., $T = 3, m = 3$), the rate is $O\left(N_n^{-\frac{p}{10(p+1)}}\right)$, which is slow, motivating the use of FuGW-based generation over naive Monte Carlo.

## 5. A practical algorithm for FuGW-based scenario tree generation

Having established that FuGW controls the MSP value gap, we now turn to the algorithmic question: given a fixed branching topology and a reference discretization of the true process, how do we find a scenario tree that minimizes the FuGW distance? This section formulates the optimization problem, proposes a Block Coordinate Descent (BCD) algorithm, analyses its convergence, and situates it among existing tree generation methods.

### 5.1 Generalities

Assume that the any generated scenario tree can be represented through a matrix $\mathrm{S} \in \mathbb{R}^{N \times Tm}$ associated to the tree's structure (where $N$ is the number of leaves, $T$ the number of stages, and $m$ the dimensionality of the decision vector at each stage) and a (probability) vector $w \in \Delta_{N-1}$, with $\Delta_{N-1}$ denoting the probability simplex ( $(N-1)$-dimensional ). Remarkably, every row of S is a branch of the tree. Consequently, our aim is to generated a scenario tree $\mathbb{T}_{\phi}(\mathrm{S}, w)$, parametrized by a matrix S and a vector $w$, which is the solution of the following optimization problem:

$$\min_{\mathrm{S} \in \mathbb{R}^{N \times Tm},\ w \in \Delta_{\mathrm{N}-1}} Fu\mathcal{GW}_{p,\alpha}\left(\overline{\mathbb{T}}_{\phi}, \mathbb{T}_{\phi}(\mathrm{S}, w)\right)$$

Although all equations are presented with respect to continuous probability distributions, they also apply to the discrete case – matrix S and vector $w$ – by simply turning integrals into summations.

**Remark 5.1: Non-convexity**. Minimizing $Fu\mathcal{GW}_{p,\alpha}$ with respect to a scenario tree parametrized through S and $w$ is a jointly non-convex problem because the structural GW term contains quadratic dependence on the coupling through the pairwise ultrametric distortion. This is the same non-convexity that appears in standard GW and QAP-type formulations, but here it is tempered by the fixed tree topology and the entropic coupling regularization used in the proposed algorithm.

**Remark 5.2: Fixed topology.** We keep the branching topology fixed and consider it as a hyperparameter. Topology selection can be handled externally by evaluating several candidates and selecting the one with the smallest final FuGW value.

### 5.2 Algorithm

We exploit the bilinear structure of the FuGW's objective by iterating the following three steps:

- computing the optimal coupling $\pi^{(k+1)}$ for the current generated tree, parametrized by the matrix $S^{(k)}$ and the probability vector $w^{(k)}$
- computing the new matrix $S^{(k+1)}$ depending on $S^{(k)}$ and $\pi^{(k+1)}$
- computing the new probability vector $w^{(k+1)}$ depending on $w^{(k)}$ and $\pi^{(k+1)}$.

---

**Algorithm 1. FGW-based Tree Generation through Block Coordinate Descent (BCD) method**

---

INPUT: a real stochastic process represented by a parametrized tree $\overline{\mathbb{T}}_{\phi}$

1. Initialize $\left(\mathrm{S}^{(0)}, w^{(0)}\right)$ randomly, $k \leftarrow 0$
2. While $k \leq$ maxIters :
   - obtain $\pi^{(k+1)} \in \underset{\pi\in\Pi(\mu_{\mathcal{A}},\mu_{\mathcal{B}})}{\operatorname{argmin}} Fu\mathcal{GW}_{p,\alpha}\left(\overline{\mathbb{T}}_{\phi}(\bar{\mathrm{S}}, \bar{w}), \mathbb{T}_{\phi}\left(\mathrm{S}^{(k)}, w^{(k)}\right)\right)$
   - compute the new $\mathrm{S}^{(k+1)}$ such that $\forall j = 1{:}\,\bar{N},\ \ \mathrm{S}_{j,*}^{(k+1)} \leftarrow \frac{\sum_{i=1}^{\bar{N}} \pi_{ij}^{(k+1)}\, \bar{\mathrm{S}}_{j,*}}{\sum_{i=1}^{\bar{N}} \pi_{i,j}^{(k+1)}}$,
   - compute the new $w^{(k+1)}$ such that $w_j^{(k+1)} \leftarrow \sum_{i=1}^{\bar{N}} \pi_{i,j}^{(k+1)}$,
   - $k \leftarrow k + 1$
3. Return the final tree $\mathbb{T}_{\phi}^{(\mathrm{maxIters})}$ as parametrized by $\mathrm{S}^{(\mathrm{maxIters})}$ and $w^{(\mathrm{maxIters})}$.

---

#### 5.3 Convergence analysis

**Monotone decrease of FuGW.** Under the assumption that $\pi^{(k+1)}$ is the global optimum of the problem $\underset{\pi\in\Pi(\mu_{\mathcal{A}},\mu_{\mathcal{B}})}{\operatorname{argmin}} Fu\mathcal{GW}_{p,\alpha}\left(\overline{\mathbb{T}}_{\phi}(\bar{\mathrm{S}}, \bar{w}), \mathbb{T}_{\phi}\left(\mathrm{S}^{(k)}, w^{(k)}\right)\right)$, then the sequence of trees generated by Algorithm 1 satisfies

$$Fu\mathcal{WG}_{p,\alpha}\left(\overline{\mathbb{T}}_{\phi}, \mathbb{T}_{\phi}\left(S^{(k+1)}, w^{(k+1)}\right)\right) \leq Fu\mathcal{WG}_{p,\alpha}\left(\overline{\mathbb{T}}_{\phi}, \mathbb{T}_{\phi}\left(S^{(k)}, w^{(k)}\right)\right) \tag{32}$$

**Fixed-point condition.** A pair $(\mathrm{S}^*, w^*)$ is a fixed point of Algorithm 1 if and only if:

- the coupling $\pi^*$ is optimal for the cost induced by $(\mathrm{S}^*, w^*)$,
- each node value is a barycenter of its assigned reference mass,
- the probabilities equal the coupling marginal.

**Remark 5.3 Comparison with nested distance based tree generation.** Exact nested-distance minimization is computationally prohibitive even for moderate stage lengths. $FuGW$ replaces it with a tractable relaxation that still controls the MSP value gap via the stability theorem. Moreover, BCD cannot be exploited with nested distance.

**6. Experiments and results**

In this paper, experiments are devoted to present the benefits provided by FuGW especially for medium-large MSP problems, specifically: T=5 with branching factor equal to 3 (leading to $3^5 = 243$ scanarios) and T=10 with branching factor equal to 2 (leading to $2^{10} = 1024$ scenarios), which could be impractical for nested distance and adapted-Wasserstein distance based approaches. Although authors are aware about methods specifically developed for dealing with a large number of scenarios and with large T values – such as those based on sampling instead of using a tree – they are not considered in this study because they belong to a different methodological framework. Here we want to show the power of FuGW in generating scenario trees whose sizes would be impractical for methods based on other distances. Anyway, some small-sized examples are also considered, such as T=2 and T=3, both with branching factor equal to 2 and 3.

The case study considered is a simple inventory management problem, with unitary cost equal to 1€. The proposed algorithm is iterated for 20 steps and restarted 5 times, from as many randomly initialized trees. Different values of $\alpha \in [0, 1]$ are analyzed, while the regularization parameter is fixed to $\varepsilon = 0.001$.

**Table 1** reports the most relevant results for the smallest case, that is T=2 and branching factor=2. For this case, the optimal value obtained on the true tree $\overline{\mathbb{T}}_\phi$ (obtained by solving the deterministic equivalent problem) is $v^*_{\mathbb{T}_\phi} = 88$ with a sum of slack values equal to 73. Columns in the table are: ($i$) the lowest FuGW value, ($ii$) the correspondent nested sinkhorn divergence (NSD) – that is a computationally efficient Sinkhorn-regularized implementation of the nested distance, ($iii$) the optimal value obtained by solving the deterministic equivalent on the generated tree $\mathbb{T}_\phi$ (specifically, the one associated to the lowest FuGW), ($iv$) the associated sum of slack values, and ($v$) the sum of slack values obtained by applying, to the actual tree, the optimal policy $x^*_{\mathbb{T}}$ obtained from the generated one (the objective value does not change). It is possible to notice that FuGW is always lower that NSD and decreases with $\alpha$ increasing. The objective value, as well as the slack values, slightly increase only for $\alpha = 1$, while NSD is not able to capture any difference between the actual tree and those generated for different values of $\alpha$; indeed, its value is always 15.5.

**Table 1.** Results for T=2 and branching factor=2. Columns are: ($i$) the lowest FuGW, ($ii$) the correspondent nested sinkhorn divergence (NSD), ($iii$) the optimal value of the deterministic equivalent for the generated tree with the lowest FuGW, ($iv$) the associated sum of slacks, and ($v$) the sum of slacks by applying, to the actual tree, the optimal policy $x^*_{\mathbb{T}}$ from the generated tree.

| $\alpha$ | FuGW | NSD | $v^*_{\mathbb{T}}$ | Slack | Slack on actual |
|---|---|---|---|---|---|
| 0.0 | 3.6085 | 15.5 | 88 | 49 | 65 |
| 0.1 | 3.4305 | 15.5 | 88 | 49 | 65 |
| 0.2 | 3.2428 | 15.5 | 88 | 49 | 65 |
| 0.3 | 3.0435 | 15.5 | 88 | 49 | 65 |
| 0.4 | 2.8303 | 15.5 | 88 | 49 | 65 |
| 0.5 | 2.5996 | 15.5 | 88 | 49 | 65 |
| 0.6 | 2.3463 | 15.5 | 88 | 49 | 65 |
| 0.7 | 2.0619 | 15.5 | 88 | 49 | 65 |
| 0.8 | 1.7320 | 15.5 | 88 | 49 | 65 |
| 0.9 | 1.3000 | 15.5 | 88 | 49 | 65 |
| 1.0 | 0.4494 | 15.5 | 89 | 56 | 69 |

**Figure 1** shows (on the left) how FuGW quickly converges to a minimum value over iterations and (on the right) how also the associated objective value decreases becoming close to that obtained for the actual tree. For a better visualization, only a subsample of the $\alpha$ values are reported in the charts.

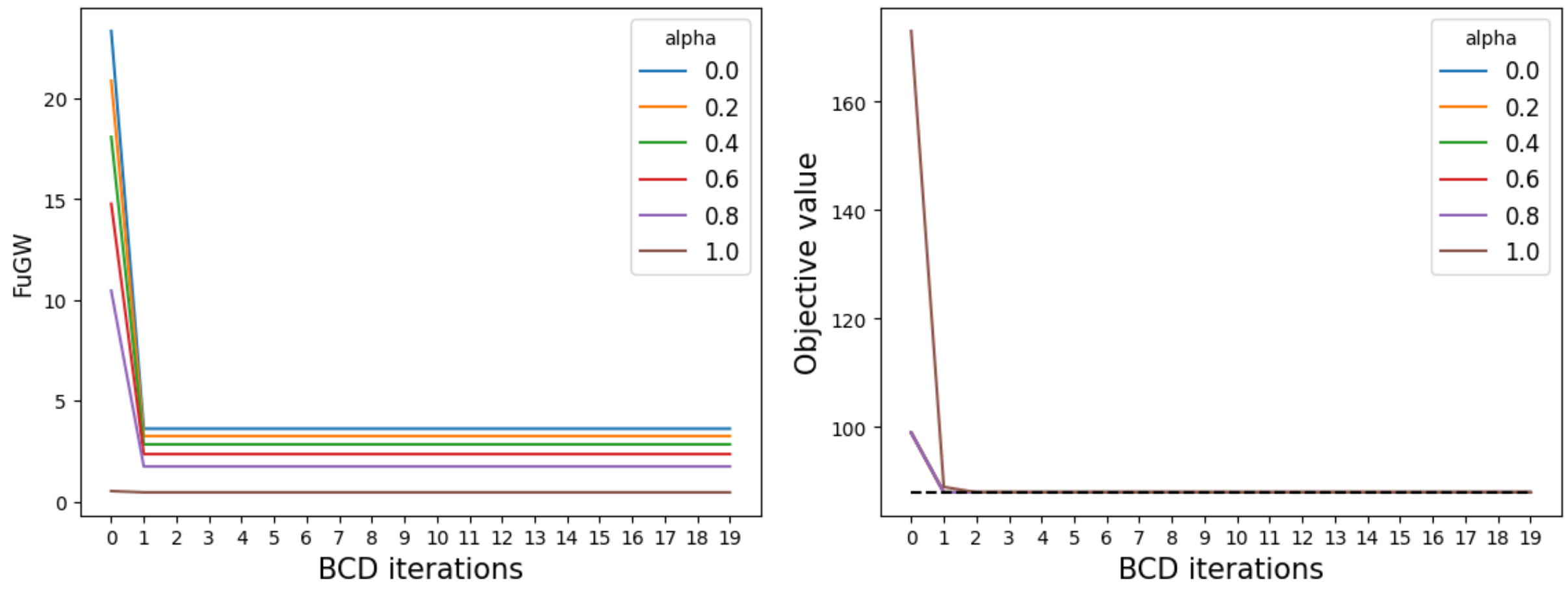

**Figure 1.** Convergence of the FuGW value (on the left) and objective value (on the right) against iterations of the proposed algorithm.

All the other case studies will present the same behaviour reported in Figure 1, so we omitted them in the next.

**Table 2** reports the most relevant results for the case with T=2 and branching factor=3. For this case, the optimal value is $v^*_{\mathbb{T}_\phi} = 113$ with a sum of slack values equal to 362. Increasing only the branching factor from 2 to 3 does not lead to different considerations with respect to the previous case; likewise for the case T=3 with branching factor=2 reported in **Table 3** (with $v^*_{\mathbb{T}_\phi} = 186$ and a sum of slack values equal to 440) and the case T=3 with branching factor=3 reported in **Table 4** (with $v^*_{\mathbb{T}_\phi} = 266$ and a sum of slack values equal to 2907).

**Table 2.** Results for T=2 and branching factor=3. Columns as described in Figure 1's caption.

| $\alpha$ | FuGW | NSD | $v^*_{\mathbb{T}}$ | Slack | Slack on actual |
|---|---|---|---|---|---|
| 0.0 | 3.6434 | 20.5 | 113 | 315 | 354 |
| 0.1 | 3.4614 | 20.5 | 113 | 315 | 354 |
| 0.2 | 3.2691 | 20.5 | 113 | 315 | 354 |
| 0.3 | 3.0649 | 20.5 | 113 | 315 | 354 |
| 0.4 | 2.8461 | 20.5 | 113 | 315 | 354 |
| 0.5 | 2.6089 | 20.5 | 113 | 315 | 354 |
| 0.6 | 2.3479 | 20.5 | 113 | 315 | 354 |
| 0.7 | 2.0540 | 20.5 | 113 | 315 | 354 |
| 0.8 | 1.7104 | 20.5 | 113 | 315 | 354 |
| 0.9 | 1.2692 | 15.5 | 112 | 279 | 345 |
| 1.0 | 0.3540 | 15.5 | 114 | 333 | 363 |

**Table 3.** Results for T=3 and branching factor=2. Columns as described in Figure 1's caption.

| $\alpha$ | FuGW | NSD | $v_{\mathbb{T}}^*$ | Slack | Slack on actual |
|---|---|---|---|---|---|
| 0.0 | 4.1608 | 15.5 | 185 | 363 | 424 |
| 0.1 | 3.9561 | 15.5 | 185 | 363 | 424 |
| 0.2 | 3.7401 | 15.5 | 185 | 363 | 424 |
| 0.3 | 3.5109 | 15.5 | 185 | 363 | 424 |
| 0.4 | 3.2656 | 15.5 | 185 | 363 | 424 |
| 0.5 | 3.0003 | 15.5 | 185 | 363 | 424 |
| 0.6 | 2.7092 | 15.5 | 185 | 363 | 424 |
| 0.7 | 2.3828 | 15.5 | 185 | 363 | 424 |
| 0.8 | 2.0039 | 15.5 | 185 | 363 | 424 |
| 0.9 | 1.5267 | 15.5 | 185 | 365 | 424 |
| 1.0 | 0.7556 | 15.5 | 186 | 342 | 432 |

**Table 4.** Results for T=3 and branching factor=3. Columns as described in Figure 1's caption.

| $\alpha$ | FuGW | NSD | $v_{\mathbb{T}}^*$ | Slack | Slack on actual |
|---|---|---|---|---|---|
| 0.0 | 4.7482 | 0.0 | 261 | 2792 | 2764 |
| 0.1 | 4.5213 | 0.0 | 261 | 2827 | 2764 |
| 0.2 | 4.2697 | 0.0 | 261 | 2792 | 2764 |
| 0.3 | 4.0109 | 0.0 | 261 | 2792 | 2764 |
| 0.4 | 3.7407 | 0.0 | 261 | 2827 | 2764 |
| 0.5 | 3.4377 | 0.0 | 261 | 2793 | 2764 |
| 0.6 | 3.1101 | 0.0 | 261 | 2793 | 2764 |
| 0.7 | 2.7416 | 0.0 | 261 | 2792 | 2764 |
| 0.8 | 2.3309 | 0.0 | 261 | 2789 | 2764 |
| 0.9 | 1.8072 | 0.0 | 261 | 2972 | 2845 |
| 1.0 | 0.6069 | 14.0 | 263 | 3015 | 2818 |

The most interesting result is that the optimal solutions obtained by using the generated trees always lead to a lower sum of slacks when applied on the actual tree compared to the sum of slacks associated to the optimal solution directly obtained from the actual tree.

Finally, **Table 5** and **Table 6** report the results for the largest case study: T=5 with branching factor=3 ($v_{\mathbb{T}_\phi}^* = 647$ with a sum of slack values equal to 68154) and T=10 with branching factor=2 ($v_{\mathbb{T}_\phi}^* = 1480$ with a sum of slack values equal to 662533). Already from **Table 5** it is possible to notice that, even if FuGW still decreases with $\alpha$ increasing, the optimal value $v_{\mathbb{T}}^*$ fluctuates in [647, 649], regardless of $\alpha$. The reason is the non-convex and multi-extremal nature of FuGW, which is further emphasized by the high dimensionality of the problem. Indeed, $v_{\mathbb{T}}^* \geq v_{\mathbb{T}}^*$, as well as the slacks.

**Table 5.** Results for T=5 and branching factor=3. Columns as described in Figure 1's caption.

| $\alpha$ | FuGW | NSD | $v^*_{\mathbb{T}}$ | Slack | Slack on actual |
|---|---|---|---|---|---|
| 0.0 | 4.3752 | 0.0 | 647 | 64886 | 68146 |
| 0.1 | 4.0021 | 20.5 | 649 | 67748 | 68632 |
| 0.2 | 3.8318 | 20.5 | 649 | 67810 | 68632 |
| 0.3 | 3.7072 | 20.5 | 649 | 67833 | 68632 |
| 0.4 | 3.6650 | 15.5 | 648 | 67231 | 68389 |
| 0.5 | 3.3949 | 0.0 | 647 | 65154 | 68146 |
| 0.6 | 3.1057 | 0.0 | 647 | 65081 | 68146 |
| 0.7 | 2.9450 | 20.5 | 649 | 67536 | 68632 |
| 0.8 | 2.6010 | 20.5 | 649 | 67700 | 68632 |
| 0.9 | 2.2760 | 20.5 | 649 | 67585 | 68632 |
| 1.0 | 1.8586 | 15.5 | 648 | 65297 | 68389 |

Finally, **Table 6** refers to the largest case study, that is T=10 and branching factor=2. Here, just 5 iterations of the algorithm has been performed, with 5 restarts, to keep the computational time similar to the previous cases (i.e., around 5 minutes by using a Google's Colab notebook). In this case the behavior of FuGW is completely the opposite of the previous experiments: FuGW increases with $\alpha$ increasing; again the reasons are the non-convexity and the high-dimensionality implied by the case study. On the other hand, for $\alpha = 0.0$ and $\alpha = 0.8$ the objective value is significantly close to $v^*_{\mathbb{T}_\phi}$ – i.e., we obtain 1483 instead of 1480 – and a quite close sum of slacks – i.e., 665617 instead of 662533. Finally, the convergence behavior, in terms of both FuGW and objective value, against the iterations of the proposed algorithm is depicted in **Figure 2**.

**Table 6.** Results for T=10 and branching factor=2.

| $\alpha$ | FuGW | NSD | $v^*_{\mathbb{T}}$ | Slack | Slack on actual |
|---|---|---|---|---|---|
| 0.0 | 7.2966 | 13.5 | 1483 | 662388 | 665617 |
| 0.2 | 7.7182 | 15.5 | 1488 | 668507 | 670737 |
| 0.4 | 8.2494 | 16.0 | 1489 | 669506 | 671761 |
| 0.5 | 8.3223 | 16.0 | 1487 | 667604 | 669713 |
| 0.6 | 8.3531 | 15.5 | 1486 | 665168 | 668689 |
| 0.8 | 8.7412 | 15.5 | 1483 | 662350 | 665617 |
| 1.0 | 9.3512 | 24.0 | 1591 | 673263 | 776209 |

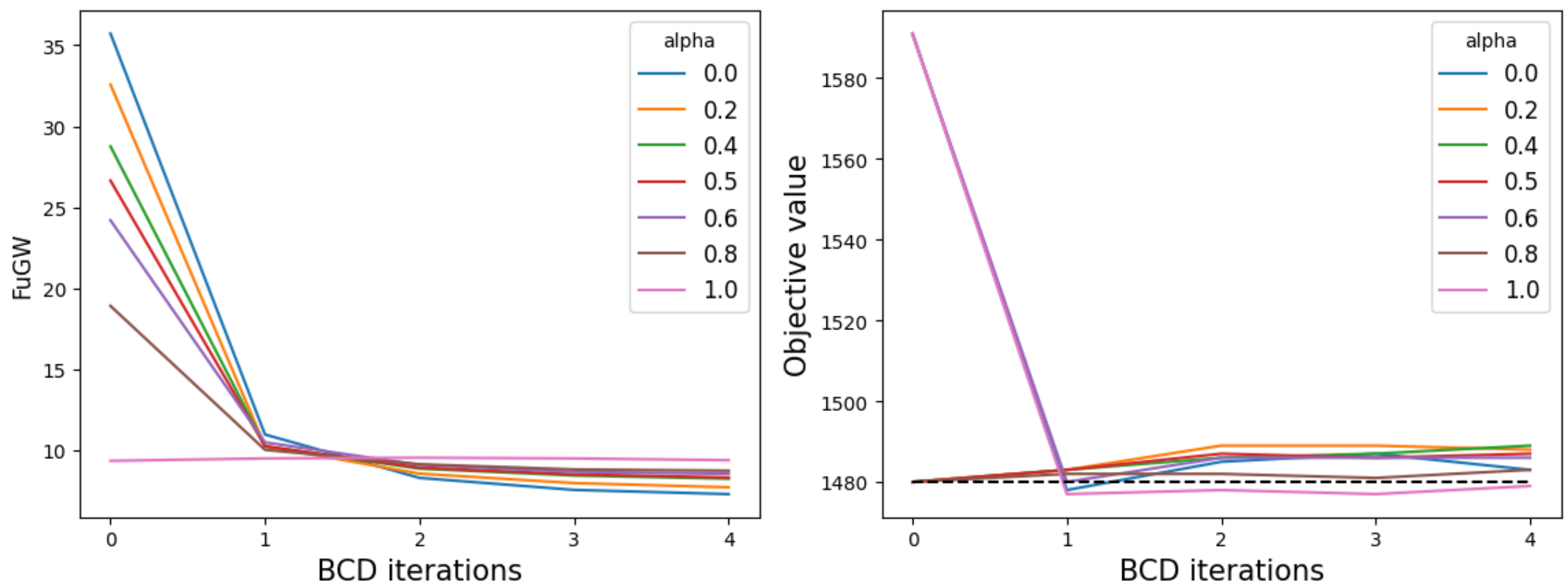

**Figure 2.** Convergence of the FuGW value (on the left) and objective value (on the right) against iterations of the proposed algorithm (case T=10 with branching factor=2).

## 7. Conclusions

We have presented a novel distances between scenario trees, namely the Fused ultrametric Gromov-Wasserstein (FuGW) which relies on the concept of ultrametric measure spaces for comparing, in a computationally efficient way, the structural properties of two scenario trees and also *fuses* it with a suitable distance between features-based representation of the scenarios. This allows us to generate a scenario tree which matches the actual one while preserving non-anticipativity. More precisely, the trade-off between fitting and violation of the non-anticipativity can be managed through the hyper-parameter $\alpha$ at the core of the proposed $Fu\mathcal{GW}_{p,\alpha}(\cdot,\cdot)$.

Empirical results on small and medium-large sized test cases of the inventory management problem prove that a simple BCD algorithm based on FuGW can quickly converge to an optimal generated tree, in a tractable computational time even on large cases (i.e., around 5 minutes on a Google's Colab notebook, for T=10, branching factor=2, given 5 BCD iterations and 5 restarts, for a single value of $\alpha$). Although the generated trees are always good, in terms of value of the objective and sum of the slacks, quality can result low when $\alpha = 1$, namely when only the structural differences between the tree are considered. Finally, the generation of the tree resulted robust – with respect to $\alpha$ – for the small-size cases, while some small fluctuations can be observed in other cases. These fluctuations are due to the strong non-convexity of the FuGW and the high-dimensionality of the associated optimization problem as implied by the size of the trees.